\documentclass[11pt]{article}

\usepackage{amssymb}
\usepackage{amsmath}
\usepackage{amsthm}
\usepackage{enumerate}
\usepackage[usenames]{color}
\usepackage{graphicx}

\usepackage[usenames]{color}
\definecolor{brown}{cmyk}{0, 0.72, 1, 0.45}
\definecolor{grey}{gray}{0.5}

\def\E{\mathbb{E}}

\title{}

\def\e{\epsilon}

   \def\p{\pi}
  \def\s{\sigma} 
 \def\om{\omega}

\newtheoremstyle{plain}%
{8pt plus2pt minus4pt}%
{8pt plus2pt minus4pt}%
{\itshape}%
{}%
{\bfseries\scshape}%
{}%
{6pt}
{}%

\newtheoremstyle{remark}%
{8pt plus2pt minus4pt}%
{8pt plus2pt minus4pt}%
{\upshape}
{}%
{\bfseries\scshape}%
{}%
{6pt}
{}%

\theoremstyle{plain}
\newtheorem{theorem}{Theorem}
\newtheorem{lemma}[theorem]{Lemma}

\theoremstyle{remark}
\newtheorem{remark}[theorem]{Remark}

\newcommand{\brac}[1]{\left(#1\right)}

\newcommand{\bfrac}[2]{\left(\frac{#1}{#2}\right)}

\def\cE{{\cal E}}

\newcommand{\beq}[1]{\begin{equation}\label{#1}}
\newcommand{\eeq}{\end{equation}}

\newcommand{\set}[1]{\left\{#1\right\}}

\def\E{\mbox{{\bf E}}}

\def\Pr{\mbox{{\bf Pr}}}
\def\whp{{\bf whp}}

\newcommand{\ignore}[1]{}

\def\e{\varepsilon}

   \def\p{\pi}
  \def\s{\sigma} 
 \def\om{\omega}

\def\E{\mbox{{\bf E}}}

\def\Pr{\mbox{{\bf Pr}}}
\def\whp{{\bf whp}}

\begin{document}
\title{Tight Hamilton Cycles in Random Uniform Hypergraphs}

\author{
{\Large{Andrzej Dudek}}\thanks{
\footnotesize {Department of Mathematics, Western Michigan University, Kalamazoo, MI 49008, \texttt{andrzej.dudek@wmich.edu}
}}
\and
{\Large{Alan Frieze}}\thanks{
\footnotesize {Department of Mathematical Sciences, Carnegie Mellon University, 
Pittsburgh, PA 15213, \texttt{alan@random.math.cmu.edu}}}~
\thanks{\footnotesize{Research partially supported by NSF grant CCF0502793}}
}

\date{\today}

\maketitle \makeatother

\begin{abstract}
In this paper we show that $e/n$ is the sharp threshold for the existence 
of tight Hamilton cycles in random $k$-uniform hypergraphs,
for all $k\ge 4$. When $k=3$ we show that $1/n$ is an asymptotic threshold. 
We also determine thresholds for the existence of other types of Hamilton cycles.
\end{abstract}

\section{Introduction}
The threshold for the existence of Hamilton cycles in the random graph $G_{n,p}$ has been known
for many years, see, e.g., \cite{AKS}, \cite{Boll} and \cite{KS}.
There have been many generalizations
of these results over the years and the problem is well understood. It is natural to try to extend
these results to hypergraphs and this has proven to be difficult. The famous P\'osa lemma \cite{P} fails to
provide any comfort and we must seek new tools.

A {\em $k$-uniform hypergraph} is a pair $(V,\cE)$ where $\cE\subseteq \binom{V}{k}$.
In the random $k$-uniform hypergraph $H_{n,p}^{(k)}$ of order~$n$ each
possible $k$-tuple appears independently with probability~$p$.

Suppose that $1\leq \ell< k$. An {\em $\ell$-overlapping Hamilton cycle} $C$
in a $k$-uniform hypergraph $H=(V,\cE)$ on $n$ vertices is a
collection of $m_\ell=n/(k-\ell)$ edges of $H$ such that for some cyclic order
of $[n]$ every edge consists of $k$ consecutive vertices and for
every pair of consecutive edges $E_{i-1},E_i$ in $C$ (in the natural
ordering of the edges) we have $|E_{i-1}\cap E_i|=\ell$. Thus, in every
$\ell$-overlapping Hamilton cycle the sets $C_i=E_i\setminus
E_{i-1},\,i=1,2,\ldots,m_\ell$, are a partition of $V$ into sets of
size $k-\ell$. Hence, $m_{\ell}=n/(k-\ell)$. We thus always assume, when discussing $\ell$-overlapping Hamilton cycles, that this necessary condition, $k-\ell$ divides $n$,  is fulfilled. In the
literature, when $\ell=k-1$ we have a {\em tight} Hamilton cycle and
when $\ell=1$ we have a {\em loose} Hamilton cycle.

A $k$-uniform hypergraph is said to be {\em $\ell$-Hamiltonian} when it contains an $\ell$-overlapping Hamilton cycle.

In a recent paper the second author proved the following:
\begin{theorem}[Frieze \cite{F}]\label{th:F}
There exists an absolute constant $c>0$ such that if $p\geq c(\log n)/n^2$ then
$$\lim_{\substack{n\to \infty\\4 |n}}\Pr(H_{n,p}^{(3)}\ \text{is $1$-Hamiltonian})=1.$$
\end{theorem}

In a subsequent paper we (essentially) extended the above theorem to $k\geq 4$.
In the following and throughout the paper,
$\omega=\omega(n)$ can be any function tending to infinity with $n$.
\begin{theorem}[Dudek and Frieze \cite{DF}]\label{th:DF}
Suppose that $k\ge 3$. If $p \ge \frac{\omega\log n}{n^{k-1}}$,
then
$$\lim_{\substack{n\to \infty\\2(k-1) |n}}\Pr(H_{n,p}^{(k)}\ \text{is $1$-Hamiltonian})=1.$$
\end{theorem}
Thus $(\log n)/ n^{k-1}$ is the asymptotic threshold for the
existence of loose Hamilton cycles, at least for $n$ a multiple of $2(k-1)$. This is because
if $p\leq (1-\e) (k-1)! (\log n) / n^{k-1}$ and $\e>0$ is
constant, then \whp\footnote{An event $\cE_n$
occurs {\em with high probability}, or \whp\ for brevity, if
$\lim_{n\rightarrow\infty}\Pr(\cE_n)=1$.} $H_{n,p}^{(k)}$ contains isolated vertices. (This follows immediately from the second moment method.)

Notice that the necessary divisibility requirement for a
$k$-uniform hypergraph to have a loose Hamilton
cycle is $(k-1) | n$. In the above two theorems we needed to
assume more, namely, $2(k-1) | n$.

As far as we know, when $\ell\geq 2$, prior to this paper, nothing
of any significance has been proven about the existence thresholds.

\begin{theorem}\label{thm:main}
\ \\[-0.3in]
\begin{enumerate}[(i)]
\item For all integers $k> \ell\geq 2$ and fixed $\e>0$, if
$p\leq (1-\e)e^{k-\ell}/n^{k-\ell}$, then
\whp\ $H_{n,p}^{(k)}$ is not $\ell$-Hamiltonian.
\item For all integers $k>\ell \ge 3$, there exists a
constant $c=c(k)$ such that if $p\geq c/n^{k-\ell}$
and $n$ is a multiple of $k-\ell$ then
 $H_{n,p}^{(k)}$ is $\ell$-Hamiltonian \whp. \label{thm:a}
\item If $k>\ell=2$ and $p\geq \omega/n^{k-2}$
and $n$ is a multiple of $k-2$,
then $H_{n,p}^{(k)}$ is $2$-Hamiltonian \whp.\label{thm:b}
\item For all fixed $\e>0$, if $k\geq 4$ and $p\geq (1+\e)e/n$, then
\whp\ $H_{n,p}^{(k)}$ is $(k-1)$-Hamiltonian, i.e. it contains a tight Hamilton cycle.
(Here $e$ is the base of natural logarithms).
\label{thm:c}
\end{enumerate}
\end{theorem}
Notice that this theorem concerns precisely those values of
$k,\ell$ not covered by Theorems \ref{th:F} and~\ref{th:DF}.

\begin{remark}
The theorem shows that $e/n$ is a sharp threshold for the existence of a tight Hamilton cycle in a $k$-uniform hypergraph, when
$k\geq 4$. When $k=3$ then the above theorem yields that  $1/n$ is an asymptotic threshold.
When $k=2$, i.e. for graphs, the sharp threshold is $\log n/n$ and as is well known, the second moment method fails.
\end{remark}

Prior to this research, we have tried combinatorial approaches
to these questions. For instance we have tried to find simple
combinatorial generalizations of Pos\'a's lemma. Surprisingly,
all it takes is the second moment method. The reason being that for $\ell>1$ the number of edges
above the density threshold is significantly larger than $n^2$. So most of weight in the mean square $\E(X^2)$ is
taken up by pairs of disjoint Hamilton cycles. For $\ell=1$ the second moment method always fails. In Section~\ref{sec:main} we present more details and give a proof of Theorem~\ref{thm:main}. Moreover, in Section~\ref{sec:pan} we consider a slightly more general problem, namely, we examine  a pancyclicity of $H_{n,p}^{(k)}$.

\begin{remark}
It is also worth mentioning that all the previous results also hold for directed random hypergraphs $\stackrel{\to}H_{n,p}^{(k)}$, where every ordered $k$-tuple appears with probability~$p$. To see this, it is enough to note that so called General Clutter Percolation Theorem of McDiarmid~\cite{MC} implies that
$$
\Pr\Big{(}\stackrel{\to}H_{n,p}^{(k)} \text{ contains an ordered Hamilton cycle}\Big{)}
\ge \Pr\Big{(}H_{n,p}^{(k)} \text{ contains a Hamilton cycle}\Big{)}.
$$
\end{remark}

\section{Proof of Theorem~\ref{thm:main}}\label{sec:main}
Let $([n],\cE)$ be a $k$-uniform hypergraph.
A permutation $\p$ of $[n]$ is an {\em $\ell$-overlapping Hamilton cycle inducing} if
$$E_\p(i)=\set{\p((i-1)(k-\ell)+j):\;j\in [k]}\in \cE\ for\ all\ i\in [n/(k-\ell)].$$
(We use the convention $\p(n+r)=\p(r)$ for $r>0$.)
Let the term {\em hamperm} refer to such a permutation.

Let $X$ be the random variable that counts the number of
hamperms $\p$ for $H_{n,p}^{(k)}$.
Every $\ell$-overlapping Hamilton cycle induces at least one hamperm and so we can
concentrate on estimating $\Pr(X>0)$.

Now
$$
\E(X) =n!p^{n/(k-\ell)}.
$$
This is because $\p$ induces an $\ell$-overlapping Hamilton cycle if and only
if a certain $n/(k-\ell)$
edges are all in $H$.

For part (i) we use Stirling's formula to argue that
$$\E(X)\leq 3\sqrt{n}\bfrac{np^{1/(k-\ell)}}{e}^n\leq 3\sqrt{n}(1-\e)^{n/(k-\ell)}=o(1).$$
This verifies part (i).

For part \eqref{thm:a} we define a constant $c=4k!ke^k$. (In order to simplify the presentation
we do not attempt to find the optimal constant.)

The beginning of the proof is the same for both part \eqref{thm:a} and \eqref{thm:b}.
We write
\begin{equation}\label{eq:expec}
\E(X)\geq \bfrac{np^{1/(k-\ell)}}{e}^n\geq c^{n/(k-\ell)}e^{-n}
\end{equation}
which goes to infinity together with~$n$.

Fix a hamperm $\p$. Let $H(\p)=(E_\p(1),E_\p(2),\ldots,E_{\p}(m_\ell))$ be the Hamilton cycle induced by
 $\p$. Then let $N(b,a)$ be the number of permutations $\p'$
such that $|E(H(\p))\cap E(H(\p^\prime))|=b$ and $E(H(\p))\cap E(H(\p^\prime))$
consists of $a$ edge disjoint paths.
Here a path is a maximal sub-sequence $F_1,F_2,\ldots,F_q$ of the edges of $H(\p)$ such that $F_i\cap F_{i+1}\neq \emptyset$
for $1\leq i<q$. The set $\bigcup_{j=1}^qF_j$ may contain other edges of $H(\p)$.
Observe that $N(b,a)$ does not depend on $\p$.

Note that
$$
\frac{\E(X^2)}{\E(X)^2} = \frac{n!N(0,0)p^{2n/(k-\ell)}}{\E(X)^2}
+ \sum_{b=1}^{n/(k-\ell)} \sum_{a=1}^{b} \frac{n!
N(b,a) p^{2n/(k-\ell)-b}} {\E(X)^2}.
$$
Since trivially, $N(0,0)\le n!$, we obtain,
\begin{equation}\label{eq:var0}
\frac{\E(X^2)}{\E(X)^2} \le 1 + \sum_{b=1}^{n/(k-\ell)} \sum_{a=1}^{b}
\frac{n! N(b,a) p^{2n/(k-\ell)-b}} {\E(X)^2}.
\end{equation}
We show that
\begin{equation}\label{eq:var}
\sum_{b=1}^{n/(k-\ell)} \sum_{a=1}^{b} \frac{n! N(b,a)
p^{2n/(k-\ell)-b}} {\E(X)^2}=
\sum_{b=1}^{n/(k-\ell)} \sum_{a=1}^{b} \frac{N(b,a)
p^{n/(k-\ell)-b}} {\E(X)}=o(1).
\end{equation}
The Chebyshev inequality implies that
$$
\Pr(X=0) \le \frac{\E(X^2)}{\E(X)^2} - 1=o(1),
$$
as required.

It remains to show \eqref{eq:var}. First we find an upper bound on $N(b,a)$.
Choose $a$ vertices $v_i$, $1\le i\le a$, on $\pi$. We have at most
\begin{equation}\label{x3}
n^{a}
\end{equation}
choices. Let
$$
b_1+b_2+\dots+b_a = b,
$$
where $b_i\ge 1$ is an integer for every $1\le i\le a$. Note that this equation has exactly
\begin{equation}\label{x2}
{b-1\choose a-1} < 2^b
\end{equation}
solutions. For every $i$, we choose a path of length $b_i$ in $H(\pi)$ which starts at $v_i$.
Suppose a path consists of edges $F_1,F_2,\ldots,F_q,\,q=b_i$. Assuming that $F_1,\dots,F_j$ are chosen, we have
at most $k$ possibilities for $F_{j+1}$. Hence, every such a path can be selected in most $k^{b_i}$ ways. Consequently, we have at most
$$
\prod_{i=1}^{a} k^{b_i} = k^b
$$
choices for all $a$ paths.

Thus, by the above consideration we can find $a$ edge disjoint paths in $H(\pi)$
with the total of $b$ edges in at most
\begin{equation}\label{eq:b}
n^a (2k)^b
\end{equation}
many ways.

Let $P_1, P_2,\dots,P_a$ be any collection of the above $a$ paths. Now we count the number of permutations $\pi^\prime$ containing these paths.

First we choose for every $P_i$ a sequence of vertices inducing this path in $\pi^\prime$. We see each edge of
$P_i$ in at most $k!$ orders. Crudely, every such sequence can be chosen in at most $(k!)^{b_i}$ ways. Thus, we have
\begin{equation}\label{eq:bb}
\prod_{i=1}^{a} (k!)^{b_i} = (k!)^b
\end{equation}
choices for all $a$ sequences.

Now we bound the number of permutations containing these sequences. First note that
$$
|V(P_i)| \ge b_i(k-\ell)+\ell.
$$
Thus we have at most
\begin{equation}\label{eq:bbb}
n-\sum_{i=1}^a (b_i(k-\ell)+\ell) = n-b(k-\ell)-a\ell
\end{equation}
vertices not in $V(P_1) \cup\dots \cup V(P_a)$. We choose a permutation
$\s$ of $V\setminus (V(P_1) \cup\dots \cup V(P_a))$. Here we have at most
\begin{equation*}
(n-b(k-\ell)-a\ell)!
\end{equation*}
choices. Now we extend $\s$ to a permutation of $[n]$. We mark
$a$ positions on $\s$ and then insert the sequences. We can do it in
\begin{equation*}
{n\choose a}a! < n^a
\end{equation*}
ways.
Consequently, the number of permutations containing $P_1, P_2,\dots,P_a$ is smaller than
\begin{equation}\label{eq:pi}
(k!)^b(n-b(k-\ell)-a\ell)! n^a.
\end{equation}
Thus, by \eqref{eq:b} and \eqref{eq:pi} and the Stirling formula we obtain
$$
N(b,a) < n^{2a} (2k!k)^b (n-b(k-\ell)-a\ell)! < n^{2a} (2k!k)^b \sqrt{2\pi n}
\left( \frac{n}{e} \right)^{n-b(k-\ell)-a\ell}(1+o(1)).
$$
Since
$$
\E(X) = n!p^{n/(k-\ell)} = \sqrt{2\pi n}
\bfrac{n}{e}^{n} p^{n/(k-\ell)} (1+o(1)),
$$
we get
$$
\frac{N(b,a) p^{n/(k-\ell)-b}} {\E(X)}
< n^{2a}  (2k!k)^b  \left( \frac{e}{n} \right)^{b(k-\ell)+a\ell}
p^{-b}(1+o(1)).$$
Finally, since $a\le b$ we estimate $e^{b(k-\ell) +a\ell} \le e^{kb}$, and consequently,
\begin{equation}\label{eq:end}
\frac{N(b,a) p^{n/(k-\ell)-b}} {\E(X)}
<\left( \frac{2k!ke^k}{n^{k-\ell}p} \right)^b \frac{1}{n^{a(\ell-2)}}(1+o(1)).
\end{equation}

Now we split the proof into two cases corresponding to two statements of Theorem~\ref{thm:main}.

\textbf{Proof of \eqref{thm:a}}:

By assumption $\ell \ge 3$ and $\frac{2k!ke^k}{n^{k-\ell}p} \le 1/2$. Thus, \eqref{eq:end} yields
$$
\frac{N(b,a) p^{n/(k-\ell)-b}} {\E(X)}
< \frac{1}{2^b n^{a}}(1+o(1))
$$
and hence
\begin{align*}
\sum_{b=1}^{n/(k-\ell)} \sum_{a=1}^{b} \frac{N(b,a) p^{n/(k-\ell)-b}} {\E(X)}
&< \sum_{b=1}^{n/(k-\ell)} \sum_{a=1}^{b} \frac{1+o(1)}{2^b n^{a}} \\
&\le \left( \sum_{b=1}^{n} \frac{1}{2^b}\right)
\left(\sum_{a=1}^n \frac{1+o(1)}{n^{a}}\right),
\end{align*}
which tends to 0 together with $n$ since the first sum is bounded by 1 and the second goes to 0.

This completes the proof of part \eqref{thm:a} of Theorem~\ref{thm:main}.

\medskip
\textbf{Proof of \eqref{thm:b}}:

Here $\ell=2$ and $n^{k-2}p\ge \omega$. Hence, we obtain in \eqref{eq:end}
$$
\frac{N(b,a) p^{n/(k-2)-b}} {\E(X)}
\le \left( \frac{2k!ke^k}{\omega} \right)^b (1+o(1)).
$$
Thus,
\beq{11x}
\sum_{b=1}^{n/(k-2)} \sum_{a=1}^{b} \frac{N(b,a) p^{n/(k-2)-b}} {\E(X)}
< \sum_{b=1}^{n} b \left( \frac{2k!ke^k}{\omega} \right)^b (1+o(1)),
\eeq
which also tends to 0 as $n$ goes to infinity, as required.

This completes the proof of part \eqref{thm:b} of Theorem~\ref{thm:main}.

\medskip
\textbf{Proof of \eqref{thm:c}}:

Let $p\geq (1+\e)e/n$. First note that as in \eqref{eq:expec} the expected value goes to infinity together with~$n$.

Next we estimate $N(b,a)$ more carefully in this case.
Suppose that $|V(P_i)|=b_i+(k-1)+t_i$. Here $t_i\ge 0$ is the number of edges of $H(\p)\setminus H(\p')$ that
are contained in $V(P_i)$.
Let $t=t_1+t_2+\cdots+t_a$. Now we argue that
\begin{equation}\label{x5}
N(b,a)\leq n^{2a}\binom{b-1}{a-1}\sum_{t\geq 0}2^{t+a}(n-b-a(k-1)-t)! (k!)^{a+t}.
\end{equation}
Here is the explanation. As before we choose $v_1,v_2,\ldots,v_a$ and $b_1,b_2,\ldots,b_a$ in
$$
n^a\binom{b-1}{a-1}
$$
ways, see \eqref{x3} and~\eqref{x2}. We then choose $t$ and then $t_1,t_2,\ldots,t_a$
in
$$
\binom{t+a-1}{a-1}<2^{t+a}
$$ ways.
Now we assign $v_i$'s to $a$ places in $\p'$ in at most
$$
n^a
$$
ways.

Now consider a fixed~$i$. We already assigned $v_i$ to a place in
$\p'$. The vertices $V(P_i)$ are now fixed but not ordered. There are
at most $k!$ ways to choose the ordering of the first edge of $P_i$.
We then choose the orderings of vertices in $V(P_i)$ induced by the $t_i$ edges
$F_1,F_2,\ldots,F_{t_i}$ in $H(\p)\setminus H(\p')$ that are contained in $V(P_i)$.
This can be done in at most $(k!)^{t_i}$ ways. Once we have ordered these edges, we claim that the ordering of any
other vertices in $V(P_i)$ are fixed. Start at the first edge, follow the ordering of
edges $F_1,F_2,\ldots,$ along
$\p'$ until we come to the first edge of $H(\p)\cap H(\p')$. The first $k-1$ of its vertices
have been ordered and so there is no choice for the $k$th vertex. Continuing in this manner
gives the claim. Therefore, we have at most
$$
\prod_{i=1}^a (k!)^{1+t_i} = (k!)^{a+t}
$$
orderings of the $V(P_i)$'s.

Having fixed the orderings of the $V(P_i)$'s and there place in $\p'$,
there are only $n-b-a(k-1)-t$ vertices left to order giving the number of choices
$$
(n-b-a(k-1)-t)!
$$
and completing the proof of \eqref{x5}.

Now we find an upper bound on every term in the summation in \eqref{x5}. Let
$$
u_t=2^{t+a}(n-b-a(k-1)-t)! (k!)^{a+t}.
$$
Then
$$
\frac{u_{t+1}}{u_t}=\frac{2k!}{n-b-a(k-1)-t}\leq \frac{1}{2}
$$
for
$$
t\leq t_0=n-b-a(k-1)-4k!.
$$
Thus,
$$
\sum_{0\le t\le t_0} u_t \le 2u_0 = 2(2k!)^a (n-b-a(k-1))!.
$$

Furthermore, for $t>t_0$ we may always assume that $t\le n-b-a(k-1)$. Hence,
\begin{align*}
\sum_{t>t_0}u_t
&\leq (2k!)^{a}(n-b-a(k-1)-t_0)! \sum_{t>t_0} (2k!)^{t}\\
&\leq (2k!)^{a}(n-b-a(k-1)-t_0)! (4k)! (2k!)^{n-b-a(k-1)}\\
&= (2k!)^a((4k)!)!(4k)!(2k!)^{n-b-a(k-1)}.
\end{align*}
But $x^m/m!\leq e^x$ for all $m\geq 0$ and so
$$
\sum_{t>t_0}u_t\leq (2k!)^a((4k)!)!(4k)!(n-b-a(k-1))! e^{2k!}.
$$
Hence,
$$\sum_{t\geq 0}(n-b-a(k-1)-t)! (2k!)^{a+t}\leq c_k(2k!)^a (n-b-a(k-1))!,
$$
where
$$c_k=2+((4k)!)!(4k)!e^{2k!}.$$
Thus,
\begin{align}
\sum_{b=1}^{n} \sum_{a=1}^{b} \frac{N(b,a) p^{n-b}} {\E(X)}
&< c_k\sum_{b=1}^{n}\frac{1}{n!p^{b}} \sum_{a=1}^{b}
n^{2a}\binom{b-1}{a-1}(2k!)^a(n-b-a(k-1))!\notag \\
&< c_k\sum_{b=1}^{n}\frac{1}{p^{b}} \sum_{a=1}^{b}
n^{2a}\binom{b-1}{a-1}(2k!)^a\bfrac{e}{n}^{b+a(k-1)}\notag \\
&=c_k\sum_{b=1}^{n}\bfrac{e}{np}^b \sum_{a=1}^{b}
\binom{b-1}{a-1}\bfrac{2k!e^{k-1}}{n^{k-3}}^a\notag \\
&=\frac{2c_kk!e^{k-1}}{n^{k-3}}
\sum_{b=1}^{n}\bfrac{e}{np}^b \brac{1+\frac{2k!e^{k-1}}{n^{k-3}}}^{b-1}\notag \\
&\leq \frac{2c_kk!e^{k-1}}{n^{k-3}}\exp\set{\frac{2k!e^{k-1}}{n^{k-4}}}
\sum_{b=1}^{n}\bfrac{e}{np}^b \label{eq:main:ddd} \\
&=o(1)\notag
\end{align}
if $np\geq e(1+\e)$.

This completes the proof of part \eqref{thm:c} of Theorem~\ref{thm:main}.

\section{Pancyclicity}\label{sec:pan}
Pancyclicity of the random graph $G_{n,p}$ has been widely studied by several researchers (see, e.g., \cite{CF,LU}). Here we consider a similar question for random hypergraphs. We say that a $k$-uniform hypergraph $H$ is {\em pancyclic} if $H$ contains a tight cycle of length~$r$,
for every $r$ in the range
$k+1\le r\le n$.

Let $X_r$ be a random variable that counts the number of tight cycles of length $r$, $k+1\le r\le n$, denoted by $C_r$. Clearly,
$$
\E(X_r) = \binom{n}{r}\frac{(r-1)!}{2} p^{r}.
$$
In particular, $\E(X_{k+1})$ goes to infinity if and only if $p\ge \frac{\omega}{n}$. We show that in general $p=\frac{\omega}{n}$ suffices for $H_{n,p}^{(k)}$ to be pancyclic.

\begin{theorem}\label{thm:pan}
\ \\[-0.3in]
\begin{enumerate}[(i)]
\item If $k=3$ and $p\geq \frac{\omega\log n}{n}$,
then $H_{n,p}^{(k)}$ is pancyclic \whp.\label{thm:pan:3}
\item If $k\ge 4$ and $p\geq \frac{\omega}{n}$,
then $H_{n,p}^{(k)}$ is pancyclic \whp.\label{thm:pan:4}
\end{enumerate}
\end{theorem}

The proof of the above theorem will be based on the following lemmas, which we prove later.
\begin{lemma}\label{lem:pan}
Let $r=o(n)$ and $p\geq \omega/n$. Then
$$
\Pr(H_{n,p}^{(k)} \text{ has no } C_r) =
\begin{cases}
o\left(1/\omega\right) & \text{for $k=3$},\\
o\left(1/(\omega n)\right) & \text{for $k\ge 4$}.
\end{cases}
$$
\end{lemma}

\begin{lemma}\label{lem:pan:2}
Let $n/\omega^{1/3} \le r\le n$ and $p\geq \omega/n$. Then
$$
\Pr(H_{n,p}^{(k)} \text{ has no } C_r) =
\begin{cases}
O\left(1/\omega^{2/3}\right) & \text{for $k=3$},\\
O\left(1/(\omega^{1/3}n)\right) & \text{for $k\ge 4$}.
\end{cases}
$$
\end{lemma}

Theorem~\ref{thm:pan} can be easily deduced from the above lemmas. Indeed, if $k\ge 4$, then the union bound yields the statement as follows.
\begin{align*}
\Pr(\exists\; k+1\le r\le n: H_{n,p}^{(k)} \text{ has no $C_r$})&\leq
\Pr(\exists\; k+1\le r\le n/\omega^{1/3}: H_{n,p}^{(k)} \text{ has no $C_r$})\\
&\quad+\Pr(\exists\; n/\omega^{1/3}\le r\le n: H_{n,p}^{(k)} \text{ has no $C_r$})\\
&= O\left(1/\omega + 1/\omega^{1/3} \right) = o(1).
\end{align*}

For $k=3$ and $p\geq \frac{\omega\log n}{n}$ we define $p'$ by $p=1-(1-p')^{\log n}$.
With this choice, we can generate $H_{n,p}^{(k)}$ as the union of $\log n$ independent copies of
$H_{n,p'}^{(k)}$. Note that $p'=\Omega(\omega/n)$. Thus, for a fixed $k+1\leq r\leq n$,
$$
\Pr(H_{n,p}^{(k)} \text{ has no $C_r$})
\leq \left( \Pr(H_{n,p'}^{(k)} \text{ has no $C_r$}) \right)^{\log n}
= O\left(1/\omega^{2/3}\right)^{\log n}
= O\left(1/\omega^{(\log n)/3}\right).
$$
Hence,
$$
\Pr(\exists\; k+1\le r\le n: H_{n,p}^{(k)} \text{ has no $C_r$})
= O\left(n/\omega^{(\log n)/3}\right) = o(1).
$$

It remains to prove Lemma~\ref{lem:pan} and \ref{lem:pan:2}.
\begin{proof}[Proof of Lemma~\ref{lem:pan}]
We make a small modification to the proof of Theorem~\ref{thm:main}. Let $C$ be a cycle of length~$r=o(n)$ and $p\geq \omega/n$. Denote by $N_r(b,a)$ the number of cycles $C'$ of length~$r$ such that $|E(C)\cap E(C')| = b$ and $E(C)\cap E(C')$ consists of $a$ edge disjoint paths. As in \eqref{eq:var0}, we obtain
$$
\frac{\E(X_r^2)}{\E(X_r)^2} \le 1 + \sum_{b=1}^{r} \sum_{a=1}^{b}
\frac{\binom{n}{r}\frac{(r-1)!}{2} N_r(b,a) p^{2r-b}} {\E(X_r)^2}
= 1 + \sum_{b=1}^{r} \sum_{a=1}^{b}
\frac{N_r(b,a) p^{r-b}} {\E(X_r)}.
$$
The Chebyshev inequality implies
\begin{equation}\label{eq:pan:2nd}
\Pr(H_{n,p}^{(k)} \text{ has no } C_r) \le \frac{\E(X_r^2)}{\E(X_r)^2} - 1 \le  \sum_{b=1}^{r} \sum_{a=1}^{b}
\frac{N_r(b,a) p^{r-b}} {\E(X_r)}.
\end{equation}

We now find an upper bound on $N_r(b,a)$. First note that as in \eqref{eq:b} we can find in $C$, $a$ edge disjoint paths with the total of $b$ edges in at most
\begin{equation}\label{eq:pan:b}
r^a (2k)^b
\end{equation}
many ways.

Let $P_1, P_2,\dots,P_a$ be any collection of the above $a$ paths, where $|E(P_i)|=b_i$ for $1\le i\le a$. Now we count the number of orderings of $V(P_1)\cup\dots\cup V(P_a)$.
As in \eqref{eq:bb} we obtain at most
\begin{equation}\label{eq:pan:bb}
\prod_{i=1}^{a} (k!)^{b_i} = (k!)^b
\end{equation}
possible orderings.

Now we bound the number of cycles $C'$ containing this collection of paths.

As in \eqref{eq:bbb} we have at most
$r-b-a(k-1)$ vertices not in $V(P_1) \cup\dots \cup V(P_a)$. We choose these vertices in at most
\begin{equation}\label{eq:pan:c}
\binom{n}{r-b-a(k-1)}
\end{equation}
many ways. Next we introduce on them a cyclic ordering having at most
\begin{equation}\label{eq:pan:cc}
(r-b-a(k-1)-1)!
\end{equation}
choices. Finally, we mark $a$ positions in order to insert paths $P_1,\dots,P_a$. Here we have at most
\begin{equation}\label{eq:pan:ccc}
(r-b-a(k-1))^a
\end{equation}
choices.

Thus, by \eqref{eq:pan:b}-\eqref{eq:pan:ccc} we get
\begin{align*}
N_r(b,a) &\le r^a (2k!k)^b n(n-1)\cdots(n-r+b+a(k-1)+1) (r-b-a(k-1))^{a-1}\\
&\le r^{2a-1} (2k!k)^b n(n-1)\cdots(n-r+b+a(k-1)+1),
\end{align*}
and consequently,
\begin{align*}
\frac{N_r(b,a) p^{r-b}} {\E(X_r)}
&\le \frac{r^{2a-1} (2k!k)^b n(n-1)\cdots(n-r+b+a(k-1)+1)p^{-b}}{n(n-1)\dots(n-r+1)/2r}\\
&\le \frac{2r^{2a} (2k!k)^b}{(n-r+b+a(k-1))\cdots(n-r+1)p^b}\\
&\le \frac{2r^{2a} (2k!k)^b}{(n-r)^{b+a(k-1)}(\omega/n)^b}\\
&= 2\left( \frac{r^2}{(n-r)^{k-1}} \right)^{a} \left( \frac{2k!k n}{(n-r)\omega} \right)^b.
\end{align*}
Hence,
\begin{align*}
\sum_{b=1}^{r} \sum_{a=1}^{b} \frac{N_r(b,a) p^{r-b}} {\E(X_r)}
&\le 2 \sum_{b=1}^{r} \sum_{a=1}^{b} \left( \frac{r^2}{(n-r)^{k-1}} \right)^{a} \left( \frac{2k!k n}{(n-r)\omega} \right)^b \\
&\le 2 \left( \sum_{b=1}^{r}  \left( \frac{2k!k n}{(n-r)\omega} \right)^b \right)
\left( \sum_{a=1}^{b} \left( \frac{r^2}{(n-r)^{k-1}} \right)^{a} \right).
\end{align*}
Since $r=o(n)$ and $k\ge 3$ both series are finite and so
\begin{align*}
\sum_{b=1}^{r} \sum_{a=1}^{b} \frac{N_r(b,a) p^{r-b}} {\E(X_r)}
= O\left(\frac{r^2}{\omega n^{k-1}}\right),
\end{align*}
which together with \eqref{eq:pan:2nd} yields  the statement of Lemma~\ref{lem:pan}.
\end{proof}

\begin{proof}[Proof of Lemma~\ref{lem:pan:2}]
Fix an $r$ such that $n/\omega^{1/3} \le r\le n$ and let $p\geq \omega/n$. Note that $H_{r,p}^{(k)}$ can be viewed as a subgraph of $H_{n,p}^{(k)}$.
Since
$$
\frac{e}{r} \le \frac{e\omega^{1/3}}{n} \ll p,
$$
the proof of Theorem~\ref{thm:main} yields that $H_{r,p}^{(3)}$ is Hamiltonian
with probability $1-O(1/\om^{2/3})$
(see \eqref{11x}) and that for $k\geq 4$, $H_{r,p}^{(k)}$ is Hamiltonian with
probability $1-O(1/(r^2p))=1-O(1/(\omega^{1/3} n))$ (see \eqref{eq:main:ddd}).
\end{proof}

\section{Concluding Remarks}

Here we summarize what is known about $\ell$-overlapping Hamilton cycles. The third column specifies the
order of magnitude of $p$ for which $H_{n,p}^{(k)}$ is $\ell$-Hamiltonian \whp.

\begin{center}
\renewcommand{\arraystretch}{2}
\begin{tabular}{| c | c | c | c | }
\hline
$\ell$ & $k$ & Order of  magnitude of $p$ & Divisibility requirement \\
\hline
\hline
$\ell=1$ & $k=3$ & $\displaystyle\frac{\log n}{n^2}$ \cite{F} & $4\ |\ n$ \\[3pt] \hline
$\ell=1$ & $k\ge 4$ & $\displaystyle\frac{\omega\log n}{n^{k-1}}$
\cite{DF} & $2(k-1)\ |\ n$ \\[3pt] \hline
$\ell=2$ & $k\ge 3$ & $\displaystyle\frac{\omega}{n^{k-2}}$ & $(k-2)\ |\ n$\\[3pt] \hline
\multicolumn{2}{|c|}{$k > \ell \ge 3$} & $\displaystyle\frac{1}{n^{k-\ell}}$ & $(k-\ell)\ |\ n$\\[3pt] \hline\hline
\multicolumn{2}{|c|}{$\ell = k-1\ge 3$}& $\displaystyle\frac{e}{n} \text{ is the sharp threshold}$ & no requirement\\[3pt] \hline
\end{tabular}
\end{center}

We close this paper with the following problems and open questions:\\[-0.3in]
\begin{enumerate}[(1)]
\item Reduce the divisibility requirement for $\ell=1$ to $(k-1)|n$.
\item Sharpen the constant for $k=3,\ell=1$ and for $k>\ell-1$.
\item Is the $\omega$ necessary in the above functions?
\item Is there a polynomial time algorithm that finds an $\ell$-overlapping Hamilton
cycle \whp\ at these densities?
\item Is $p=(\omega\log n)/n$ the threshold for pancyclicity of $H_{n,p}^{(3)}$? (Most likely the $\log n$ factor is unnecessary.)
\end{enumerate}

\section{Acknowledgment}
We would like to thank the referee for his or her valuable comments and suggestions.

\end{document}